\documentclass{article}
\newtheorem{fakt}{Lemma}[section]
\newtheorem{f}[fakt]{Fact}

\newtheorem{tw}[fakt]{Theorem}
\newtheorem{wtw}{Theorem}
\newtheorem{cor}[fakt]{Corollary}
\newtheorem{prop}[fakt]{Proposition}
\newcommand {\re}{{\it Remark. }}
\newcommand {\dow}{{\it Proof:} }
\newcommand {\kon}{\hfill $\diamond$ \\}
\newcommand {\bN}{{\bf N}}
\newcommand {\bR}{{\bf R}^n}
\newcommand {\lk}{\left\{}
\newcommand {\rk}{\right\}}
\newcommand {\mr}{\mathrm}

\usepackage{pictexwd,dcpic,epsf}

\begin {document}
\title {Connectedness at infinity of systolic complexes and groups}
\author{Damian Osajda
\footnote{Author is a fellow of the Marie Curie
Intra-European Fellowship, contract MEIF CT2005 010050 and was
partially supported by the Polish Scientific Research Commitee
(KBN) grant 2 P03A 017 25}}
\maketitle
\begin{center}
{\it Instytut Matematyczny, Uniwersytet Wroc{\l}awski

pl. Grunwaldzki 2/4, 50-384 Wroc{\l}aw, Poland

{\rm and}

Institut de Math\'ematiques de Jussieu, Universit\'e Paris 6

Case 247, 4 Place Jussieu, 75252 Paris Cedex 05, France

\vspace{0.3cm}
{\rm E-mail:} {\tt dosaj@math.uni.wroc.pl}}
\end{center}

\begin{abstract}
By studying connectedness at infinity of systolic groups we
distinguish them from some other classes of groups, in particular
from the fundamental groups of manifolds covered by euclidean
space of dimension at least three. We also study
semistability at infinity for some systolic groups.


\end{abstract}
\section{Introduction}
Systolic complexes were introduced by T. Januszkiewicz and
J. \'Swi{\c a}tkowski
(\cite{JS1}) and, independently, by F. Haglund (\cite{H}) as
combinatorial analogues of nonpositively curved spaces.
They are simply connected simplicial complexes satisfying some
local combinatorial conditions: roughly speaking there is a lower bound
for the lengths of ``essential" closed paths in the one-skeleton of every
link.

This condition is an analogue of the Gromov condition implying
nonpositive curvature
for cubical complexes. However systolic complexes
equipped with the metric for which every simplex is isometric to the
regular euclidean simplex are not necessarily nonpositively curved. Conversely,
there exist
nonpositively curved spaces (e.g. manifolds of dimension at least three) that
do not admit systolic triangulations. Nevertheless systolic
spaces (sometimes referred to as complexes of {\it simplicial nonpositive
curvature} - SNPC) possess many properties analogous to those of
spaces of non-positive curvature (we use \cite{JS1} as a reference
for this theory).
They are contractible (this is an analogue of the Cartan-Hadamard theorem),
with some additional assumptions they are Gromov
hyperbolic or $CAT(0)$ ($CAT(-1)$),
and complexes of groups with links satisfying the same conditions
as links in systolic spaces
are developable.

The latter property allows one to construct many examples of systolic
spaces and groups (i.e. groups acting geometrically on systolic
complexes) with some additional properties and to
answer some open questions.

Systolic groups are biautomatic (\cite{JS1}). They can be
considered as higher dimensional analogues of small cancellation
groups (\cite{W}). Ideal boundaries of $7$-systolic groups
are strongly hereditarily aspherical compacta (\cite{O}).

We study the topology at infinity of
systolic complexes and groups. It turns out that connectedness at
infinity allows to distinguish systolic groups from some other
classes of groups,
in particular from fundamental groups of closed manifolds covered
by ${\bf R}^n$, for $n>2$ (see Section \ref{hcis}). This extends some
of the results obtained in \cite{JS2} and \cite{W} (see Subsection \ref{nex}
below). Our main tool is the property $Conn_k^{\infty}(G)$ (see Section
\ref{sci}) saying that for every contractible,
rigid $G-CW$ complex $X$ on which $G$ acts properly
discontinuously and cocompactly and for every neighborhood $U$ of
infinity in $X$ there exists a smaller neighborhood of infinity $V\subset U$
such that every $k$-dimensional sphere in $V$ is homotopically
trivial in $U$. The main result is the following.

\begin{wtw}[Theorem \ref{tw1} in Section \ref{hcis}]
\label{A}
Let $X$ be a locally finite systolic complex, let $G$ be a systolic
group and let $k\geq
2$ be a natural number. Then the conditions $Conn_k^{\infty}(X)$
and $Conn_k^{\infty}(G)$ hold.
\end{wtw}

As a corollary of this theorem, using some homological
algebra, we obtain the following.

\begin{wtw}[Theorem \ref{tw3} in Section \ref{hcis}]
\label{AA}
Systolic groups are not simply connected at infinity.
\end{wtw}


In Subsection \ref{nex} we give some examples of groups that
are not systolic. Among them there are some Coxeter groups, lattices in
isometry groups of some buildings and nonpositively curved cubical
complexes.

In Section \ref{scis} we
study topology at infinity of locally finite systolic chamber complexes
such that links of vertices are connected and complements of
open balls of radius $2$ in those links are connected--- condition
$R(v,X)$ in Section \ref{scis}. Our results can be
viewed as analogues of some
of the results obtained in \cite{BMcCM} for $CAT(0)$
complexes
and in \cite{DM} for Coxeter groups.
We prove the following.

\begin{wtw}[Theorem \ref{tw}, Section \ref{hcis}]
\label{B}
Let $X$ be a locally finite systolic chamber complex of
dimension greater than one. Suppose that
$R(v,X)$ holds
for almost every vertex $v\in X$ (i.e. for all but finitely many).
Then $X$ is not
simply connected at infinity. Moreover, if $G$ acts properly and
cocompactly on $X$ then $G$ is semistable at infinity.
\end{wtw}

As shown in Section \ref{scis} for a systolic group
acting on a normal pseudoma\-ni\-fold $X$ the
condition $R(v,X)$ holds for every vertex $v$.
Such groups of arbitrarily large cohomological dimension are constructed
in \cite{JS1}.
Those groups are the only systolic groups of cohomological dimension
greater than two
known to us at the moment.

Semistability at infinity is conjectured to hold for all one ended
finitely presented groups. Thus we just give new examples of
groups for which it holds.

I would like to thank Tadeusz Januszkiewicz and Jacek \'Swi{\c a}tkowski for
introducing me to the subject and suggesting the problem. I am also
grateful to them and to Jan Dymara, Ross Geoghegan, Peter Kropholler,
Ian Leary and Carrie Schermetzler for helpful
conversations.

\section{Preliminaries}
\label{pre}

\subsection{Systolic complexes and groups}
\label{scg} A simplicial complex $X$ is {\it flag} if
every finite subset of its vertices pairwise connected by edges
spans a simplex of $X$. Following \cite{JS1} and \cite{JS2} we say
that, for a given natural $k\geq 4$, a simplicial complex $X$ is
{\it k-large} if it is flag and if every cycle $\gamma$ in $X$ (i.e.
its subcomplex homeomorphic to the circle) of length $4\leq
|\gamma| <k$ has a diagonal (i.e. an edge connecting two
nonconsecutive vertices in $\gamma$). Here $|\gamma|$ denotes the
number of edges of $\gamma$.

A simplicial complex $X$ is {\it locally k-large} if for every
simplex $\sigma\neq \emptyset$ of $X$ its link $X_{\sigma}$ in $X$
is $k$-large.

$X$ is {\it k-systolic} if it is locally $k$-large, connected and
simply connected.

Because $k=6$ is of special importance in that theory,
$6$-systolic complexes are called {\it systolic}.

A group acting geometrically (i.e. properly discontinuously and
cocompactly) by simplicial automorphisms on a $k$-systolic (resp.
systolic) complex is called {\it k-systolic} (resp. {\it
systolic}).

For example free groups and 
fundamental groups of surfaces are systolic. For other
examples see below.

In the rest of this subsection we list some results concerning
systolic complexes and groups which we will use later. 
If it is not stated otherwise
they can be found in \cite{JS1} and we follow the
notation of that paper. In particular, for a simplicial complex
$X$ we denote by $X^{(i)}$ its $i$-skeleton, and we denote
by $\sigma \ast
\rho$ the join of simplices $\sigma$ and $\rho$.

Recall that a subcomplex $Y$ of a simplicial complex $X$ is {\it
full} if every set $A$ of vertices of $Y$ spanning a simplex of
$X$ spans a simplex in $Y$. The following facts are immediate:

\begin{f}
\label{f0.1}
\begin{enumerate}
\item
A full subcomplex in a (locally) $k$-large complex is (locally)
$k$-large.
\item
Links of a $k$-large complex are $k$-large.
\item
(\cite{JS1}, 1.8.5) There is no $k$-large
triangulation of the $2$-sphere for $k\geq 6$. Hence no
triangulation of a manifold of dimension $n\geq 3$ is $6$-large
since $2$-spheres would occur as links of some simplices in
such triangulation.
\end{enumerate}
\end{f}

Now let $X$ be a systolic complex and $\sigma$ its simplex. By
Section 7 of \cite{JS1} one can define a {\it closed combinatorial
ball of radius $i$ around $\sigma$ in $X$}, $B_i(\sigma,X)$,
inductively: $B_0(\sigma,X)=\sigma$ and $B_i(\sigma,X)=\bigcup \lk
\tau : \tau \cap B_{i-1}(\sigma,X)\neq \emptyset \rk$, for $i>0$.
(Note that in Section 7 of \cite{JS1} the combinatorial balls and
the deformation retractions we consider are defined more
generally, for any {\it convex} subcomplex of $X$.)

By $S_i(\sigma,X)$ we denote the subcomplex of
$B_i(\sigma,X)$ spanned by the vertices at combinatorial distance
$i$ from $\sigma$, i.e. not belonging to $B_{i-1}(\sigma,X)$.
By ${\stackrel{\circ}{B_i}(\sigma,X)}$ we denote the {\it interior} of the closed
combinatorial $i$-ball around $\sigma$ in $X$, i.e.
${\stackrel{\circ}{B_i}(\sigma,X)}=B_i(\sigma,X)\setminus
S_i(\sigma,X)$.

Thus one can define closed combinatorial balls of small radii in
$k$-large complexes and they are isomorphic to ones in the
corresponding universal covers. In the following lemma the first
claim implies the second one, in view of Fact \ref{f0.1}.

Proofs of the following two lemmas (except the last part of the 
second lemma)
can be found in Section 7 of 
\cite{JS1}.

\begin{fakt}
\label{l0.2}
The ball $B_i(\sigma,X)$ and the sphere $S_i(\sigma,X)$ are
full subcomplexes of $X$ and they are $k$-large.
\end{fakt}

\begin{fakt}
\label{l0.3} For every simplex $\tau
\in
S_i(\sigma,X)$, $i>0$, $\rho = S_{i-1}(\sigma,X)\cap X_{\tau}$ is a
single simplex and $X_{\tau}
\cap B_{i}(\sigma,X)= B_1(\rho, X_{\tau})$
and $X_{\tau}
\cap S_{i}(\sigma,X)= S_1(\rho, X_{\tau})$.
\end{fakt}

\dow
We prove here only that 
$X_{\tau}
\cap S_{i}(\sigma,X)= S_1(\rho, X_{\tau})$.

Let a vertex $v$ belong to $X_{\tau}
\cap S_{i}(\sigma,X)$. Then, since $\rho \subset
B_{i-1}(\sigma,X)$ we have $v\in S_1(\rho, X_{\tau})$ and
hence $X_{\tau}
\cap S_{i}(\sigma,X)\subset S_1(\rho, X_{\tau})$.

Now let a vertex $v$ belong to $S_1(\rho, X_{\tau})$.
Then it is obvious that $v\in S_{i}(\sigma,X) \cup
S_{i-1}(\sigma,X)$. Assume $v\in
S_{i-1}(\sigma,X)$. Then $v\ast \rho \in
S_{i-1}(\sigma,X)\cap X_{\tau}=\rho$ which is a contradiction.
Hence $X_{\tau}
\cap S_{i}(\sigma,X)\supset S_1(\rho, X_{\tau})$.
\kon

In the rest of the paper we call the simplex $\rho$, as in the
above lemma the
{\it projection of $\tau$ on $B_{i-1}(\sigma,X)$}.

For a simplicial complex $X$, we denote by $X'$ its first barycentric
subdivision. For a simplex $\sigma \in X$ we denote by $b_{\sigma}$ the
barycenter of $\sigma$.

By the lemma above we can define an {\it elementary contraction}
$$
\pi_{B_{i}(\sigma,X)}:B_{i+1}(\sigma,X)'\to B_{i}(\sigma,X)'
$$
between barycentric subdivisions of balls by putting
$$
\pi_{B_{i}(\sigma,X)}(b_{\nu})=\lk \begin{array}{ccc} b_{\nu \cap
B_{i}(\sigma,X)} & if & \nu \cap B_{i}(\sigma,X)\neq \emptyset \\
b_{X_{\nu} \cap B_{i}(\sigma,X)} & if & \nu \cap B_{i}(\sigma,X)=
\emptyset
\end{array}
\right .
$$
and then extending simplicially. In Section 8 of \cite{JS1} it is
shown that $\pi_{B_{i}(\sigma,X)}$ is a deformation retraction
and that $ \pi_{B_{i}(\sigma,X)}(B_{i+1}(\sigma,X)\setminus
\stackrel{\circ}{B_{i}}(\sigma,X))\subset
S_{i}(\sigma,X)$.

Then we define a deformation retraction
$P_{B_{i}(\sigma,X)}:X\to B_{i}(\sigma,X)$ as follows. If $x\in
B_j(\sigma,X)$ then
$P_{B_{i}(\sigma,X)}(x)=\pi_{B_{i}(\sigma,X)}\circ
\pi_{B_{i+1}(\sigma,X)}\circ ... \circ
\pi_{B_{j-1}(\sigma,X)}(x)$.

\begin{fakt}
\label{l0.4} For $j>i$, the projection
$P_{B_{i}(\sigma,X)}|_{B_{j}(\sigma,X)}:{B_{j}(\sigma,X)}X\to
B_{i}(\sigma,X)$ provides a deformation retraction of
$B_{j}(\sigma,X)\setminus {\stackrel{\circ}{B_i}(\sigma,X)}$ onto
$S_i(\sigma,X)$ within $B_{j}(\sigma,X)\setminus
{\stackrel{\circ}{B_i}(\sigma,X)}$.
\end{fakt}

For $k\geq 6$, by Corollary 1.5 of \cite{JS1} (compare Fact
\ref{f0.1} above), a simplicial complex $Y$ is $k$-large iff it is
locally $k$-large and the minimal length of a homotopically
nontrivial (in $Y$) loop in the $1$-skeleton $Y^{(1)}$ of $Y$ is at
least $k$. Hence we obtain the following.

\begin{fakt}
\label{l0.15} Let $k\geq 6$ and let $Y$ be a $k$-large simplicial
complex. If $p:X\to Y$ is the universal cover of $Y$ and
$m<\frac{k-1}{2}$ then for $i=0,1,2,...,m$ the map
$p|_{B_i(\sigma,X)}:B_i(\sigma,X)\to p(B_i(\sigma,X))$ is an
isomorphism.
\end{fakt}

By part 3 of Fact \ref{f0.1} there is no systolic manifold above
dimension two. But there are systolic groups acting geometrically
on spaces that
are "close" to manifolds as we describe below.

A simplicial complex $X$ is called a {\it chamber complex of
dimension $n$} if it is the union of its $n$-simplices (which are
called {\it chambers} of $X$) and for every $(n-1)$-dimensional
face of $X$ there exist at least two chambers containing it. It is
easy to see that links in a chamber complex are themselves chamber
complexes. A {\it gallery} in a chamber complex is a finite
sequence of maximal simplices such that two consecutive simplices
share a common face of codimension $1$. A chamber complex is said
to be {\it gallery connected} if any two chambers can be connected
by a gallery. A chamber complex is {\it normal} if it is gallery
connected and all its links of dimension above $0$ are gallery
connected. A chamber complex is a {\it pseudomanifold} if every
codimension one face belongs to exactly two maximal simplices.

In \cite{JS1} it is shown (Corollary 19.2) that for every natural
$n$ and every $k\geq 6$ there exists an $n$-dimensional compact
chamber complex that is $k$-large. Such examples are constructed
as developments of some simplices of groups. Moreover, they can
be pseudomanifolds. In Section \ref{extra} we show that the groups
constructed this way are one--ended and semistable at infinity.
Here we give some background
on the constructions given in \cite{JS1}.

Recall (we use here the notation of \cite{JS1}; a standard
reference for complexes of groups is \cite{BH}) that for a given
simplex $\Delta$ a {\it simplex of groups ${\cal G}$ over
$\Delta$} is a family of groups
$\lk G_{\sigma}|\; \sigma \; \;{\mr {is \;\; a \;\; subsimplex \;\;
of}}\;\; \Delta \rk$ together with the family of {\it structure
homomorphisms} $\lk \psi_{\sigma \tau}:G_{\sigma}\to G_{\tau}|\;
\tau \subset \sigma \rk$ satisfying natural commutation
relations. A simplex of groups ${\cal G}$ is called {\it
developable} if there exists a simplicial complex $X$ ({\it
development}) and a group $G$ acting on $X$ such that $G\backslash
X=\Delta$ where $\Delta$ is identified with some maximal simplex
of $X$, \\ $G_{\sigma}=Stab(\sigma,G)$ for $\sigma$ being a
subsimplex of $\Delta$ and $\psi_{\sigma \tau}$ are the
inclusions $G_{\sigma} \to G_{\tau}$ for $\tau \subset \sigma$.
We allow the case $\sigma=\emptyset$ and then we put $G_{\emptyset}=
G$.
Then we write ${\cal G}=G\backslash \backslash X$.
We say also that $X$ is the development of ${\cal G}$ with respect
to the morphism $m:{\cal G}\to G$. The morphism $m$ is the family
of injections $m_{\sigma}:G_{\sigma}\to G$ agreeing with structure
homomorphisms, for $\sigma$ being simplices
of $\Delta$.
Moreover, if
${\cal G}$ is developable then there exists a simply connected
simplicial complex $\widetilde {\cal G}$ ({\it universal development
of ${\cal G}$}) and a group $\widehat {\cal G}$ ({\it direct limit of
${\cal G}$}) such that ${\cal G}=\widehat {\cal G}\backslash
\backslash \widetilde {\cal G}$. For a subsimplex $\sigma$ of $\Delta$
one can consider a simplex of groups ${\cal G}^{\sigma}$ over a
link $\Delta_{\sigma}$, being a restriction of ${\cal G}$ to the
link". There is a canonical morphism $i_{\sigma}:{\cal G}^{\sigma}
\to G_{\sigma}$.

Following \cite{JS1} for a natural $k\geq 4$ we will call a
simplex of groups ${\cal G}$, {\it locally $k$-large} if for every
link $\Delta_{\sigma}$, the development of ${\cal G}^{\sigma}$ with respect
to the canonical morphism $i_{\sigma}:{\cal G}^{\sigma}
\to G_{\sigma}$ is $k$-large.

\begin{prop}[Corollary 17.4 in \cite{JS1}]
\label{p0.5} For $k\geq 6$ every locally $k$-large simplex of
groups is developable.
\end{prop}

In Section 18 of \cite{JS1} a notion of {\it (local)
extra-tilability} of simplices of groups and of their
developments is introduced. The most important for us is the
following property.

\begin{prop}[Proposition 18.3 in \cite{JS1}]
\label{p0.6} Let ${\cal G}$ be a locally extra-tilable simplex of
groups. Then the action of the direct limit $\widehat {\cal G}$ on the
universal development $\widetilde {\cal G}$ has the following
property: each $i$-ball $B_i(\sigma,\widetilde {\cal G})$, for any
natural $i$, is a strict fundamental domain for the action of
a unique
subgroup of $\widehat {\cal G}$.
\end{prop}

By applying Proposition \ref{p0.6}, it is shown in \cite{JS1}
that for given natural $k\geq 6$, simplex $\Delta$ and a family of
finite groups $ \lk A_s|\; s\; \; {\mr {is \; \; a \; \;
codimension \;\; 1 \;\; face \;\; of }}\;\; \Delta \rk $
there exists
a
simplex of groups ${\cal G}$ over $\Delta$ such that $G_s=A_s$,
for every codimension one face, and ${\cal G}$ admits finite $k$-large
development. The explicit
construction of local groups $G_{\sigma}$ of such ${\cal G}$
goes by induction on the
codimension of $\sigma$. The important point is that for every
$\sigma$, local developments of ${\cal G}$ at $\sigma$ are $k$-large and
extra-tilable.

\subsection{Homology of groups}
The following two propositions and their proofs
were supplied by Ian Leary.
The main reference for homology of groups is \cite{B2}.

\begin{prop}
\label{I1}
If a group $G$ acts with finite stabilizers on a finite dimensional
contractible simplicial complex $X$ then:
\begin{enumerate}
\item ${\mr {cd}}_{\bf Q}G\leq {\mr{dim}}(X)$;
\item There exists a free ${\bf Q}G$-module $F=\bigoplus _i{\bf Q}G$ and
$n'\leq {\mr {dim}}(X)$ such, that $H^{n'}(G;F)\neq 0$;
\item There exists a free ${\bf Z}G$-module $\widetilde{F}=
\bigoplus_j{\bf Z}G$ and $n=n'$ or $n=n'+1$ such, that
$H^{n}(G;\widetilde F)\neq 0$.
\end{enumerate}
\end{prop}
\dow
{\sl 1.} Let $H$ be a finite subgroup of $G$ and let
$\pi:{\bf Q}G\to {\bf Q}G/H$ be the projection map given
by $\pi(g)=gH$. Then there is a section
$s:{\bf Q}G/H\to {\bf Q}G$ of $\pi$, defined by
$s(gH)=\frac{1}{|H|}\sum_{h\in H}gh$. Hence ${\bf Q}[G/H]$
is ${\bf Q}G$-projective and the simplicial chain complex
$C_{\ast}(X;{\bf Q})$ is a projective resolution of ${\bf Q}$
over ${\bf Q}G$ (compare Exercise I.8.5 in \cite{B2}).

{\sl 2.} Since ${\mr {cd}}_{\bf Q}G\leq{\mr {dim}}(X)<\infty$,
we have
$$
{\mr {cd}}_{\bf Q}G={\mr {max}}\lk n'\;|\; \exists \;{\mr
{a}} \; {\bf Q}G-{\mr{module\;}} M {\mr{with}} \; H^{n'}(G;M)\neq 0 \rk.
$$
Given such an $M$, let $F$ be a free ${\bf Q}G$-module
fitting into the short exact sequence
$$
0\to K\to F\to M\to 0.
$$
The cohomology long exact sequence gives
$$
H^{n'}(G;F)\to H^{n'}(G;M)\to H^{n'+1}(G;K).
$$
Since $H^{n'+1}(G;K)=0$ we have $H^{n'}(G;F)\neq 0$.

{\sl 3.} There is a short exact sequence of abelian groups
$$
0\to \bigoplus_{m\geq 1}{\bf Z}\to \bigoplus_{m\geq 1}{\bf Z}\to
{\bf Q}\to 0.
$$
Hence there is a short exact sequence of ${\bf Z}G$-modules
$$
0\to \bigoplus_{m\geq 1}{\bf Z}G\to \bigoplus_{m\geq 1}{\bf Z}G\to
{\bf Q}G\to 0
$$
and there is a short exact sequence of ${\bf Z}G$-modules
$$
0\to \widetilde{F}\to \widetilde{F}\to F\to 0.
$$
The cohomology long exact sequence gives
$$
H^{n'}(G;\widetilde{F})\to H^{n'}(G;\widetilde{F})\to H^{n'}(G;F)\to
H^{n'+1}(G;\widetilde{F}).
$$
Since $H^{n'}(G;F)\neq 0$ we obtain $H^{n'}(G;\widetilde{F})\neq 0$
or $H^{n'+1}(G;\widetilde{F})\neq 0$.
\kon

\begin{prop}
\label{I2}
If a group $G$ acts with finite stabilizers and cocompactly on a
finite dimensional contractible simplicial complex $X$ then
there exists a natural number $n\leq {\mr {dim}}(X)$ such,
that $H^{n}(G;{\bf Z}G)\neq 0$.
\end{prop}
\dow
Let $C_{\ast}(X)$ denote a simplicial chain complex for $X$ and
let $P_{\ast}$ be a free resolution for ${\bf Z}$ over
${\bf Z}G$. Consider the spectral sequence for
$$
E^{ij}_0={\mr {Hom}}_G(C_i(X)\otimes P_j,M),
$$
where $M\cong\bigoplus_k {\bf Z}G$ (compare Chapter VII.7 in
\cite{B2}).
For each fixed $i$, $C_i(X)\otimes P_{\ast}$ is a free resolution
of $C_i(x)$, so in the spectral sequence in which the differential
$d_0:E^{ij}_0\to E^{i,j+1}_0$ is induced by the differential
on $P_{\ast}$, we have $E^{ij}_1\cong {\mr {Ext}}_G^j(C_i,M)$.
By Shapiro's lemma (cf. e.g. Proposition III.6.2 in \cite{B2}),
${\mr {Ext}}_G^{\ast}({\bf Z}G/H,M)\cong {\mr {Ext}}_H^{\ast}
({\bf Z},M)=H^{\ast}(H;M)$, for a subgroup $H<G$.
By the hypotheses, for a stabilizer $H$ of every
simplex in $X$, we have $H^j(H;M)=0$, for $j>0$. Hence
${\mr {Ext}}_G^j(C_i,M)=0$, for $j>0$ and we have
$E^{ij}_1=0$, for $j>0$, and $E^{i0}_1={\mr {Hom}}_G
(C_i,M)$. Thus the spectral sequence collapses at
$E_2$, with $E^{ij}_2=0$, for $i>0$ and $E^{i0}_2=H^i_G(X;M)$.
On the other hand, the total complex for $C_{\ast}\oplus
P_{\ast}$ is a free resolution for ${\bf Z}$ over
${\bf Z}G$, and so the spectral sequence must converge to
a filtration of $H^{\ast}(G;M)$.
Thus we have a natural
isomorphism between $H^{\ast}{\mr {Hom}}_G(C_{\ast}(X);M)$
and $H^{\ast}(G;M)$ (compare Exercise VIII.7.4 in \cite{B2}).

Since $C_{\ast}(X)$ is a chain complex of finitely generated
$G$-modules we have
$$
{\mr {Hom}}_G(C_{\ast}(X),\bigoplus_j
{\bf Z}G)=\bigoplus_j {\mr {Hom}}_G(C_{\ast}(X),
{\bf Z}G).
$$
Now taking cohomology we obtain
$$
H^{\ast}(G;M)=H^{\ast}{\mr {Hom}}_G(C_{\ast}(X);M)=
$$
$$
=\bigoplus_k H^{\ast}{\mr {Hom}}_G(C_{\ast}(X);{\bf Z}G)=
\bigoplus_k H^{\ast}(G;{\bf Z}G).
$$

Using Proposition \ref{I1} and setting $M=\widetilde F$ we have that
there exists $n$ such that
$H^n(G;{\bf Z}G)\neq 0$.
\kon

\subsection{Connectedness and acyclity at infinity}
\label{sci} For a topological space $Y$ and an integer $n\geq -1$
denote by $Conn_n^{\infty} (Y)$ the following condition:

{\it For every compact
$K\subset Y$ there exists a compact $L\subset Y$ such that
$K\subset L$ and every map $S^n=\partial B^{n+1} \to Y\setminus
L$ extends to a map $B^{n+1}\to Y\setminus K$.}
This condition is also called {\it vanishing of the $n$-th homotopy
pro-group at infinity}---see \cite{G}.

If $Conn_k^{\infty}(Y)$ holds for every $k\leq n$ then we say that
$Y$ is {\it n-connected at infinity} (compare e.g. \cite{G}).
$Y$ is $(-1)$-connected at infinity iff it is not compact.
$Y$ is $0$-connected at infinity iff it has one end.
A space $1$-connected at infinity is also called {\it simply
connected at infinity}.

Recall that two maps $f_1,f_2:X\rightarrow Y$ are said to be {\it
properly homotopic} if there exists a proper (i.e. with compact
preimages of compact
sets) homotopy $F:X\times [0,1]\rightarrow Y$ 
joining them. A proper map $f:X\rightarrow Y$ is a {\it
proper homotopy equivalence} if there exists a proper map
$g:Y\rightarrow X$ such that $f\circ g$ and $g\circ f$ are
properly homotopic to, respectively $id_Y$ and $id_X$. It is
obvious that if $f:X\rightarrow Y$ is a proper homotopy
equivalence then $Conn_n^{\infty} (X)$ iff $Conn_n^{\infty} (Y)$.
A map $f \colon X \to Y$ between CW complexes is {\it {CW-proper}} if for each $n$ there
exists $k$ such that $f(X^{(n)})\subset X^{(k)}$ and $f|_{X^{(n)}}\colon X^{(n)}\to Y$ is proper.
A definition of {\it {CW-proper homotopy equivalence}} is then analogous to 
proper homotopy equivalence. If $f\colon X\to Y$ is a CW-proper homotopy
equivalence then $Conn_n^{\infty} (X)$ iff $Conn_n^{\infty} (Y)$.

The next theorem is a slight generalization of a part of Theorem
5.8.2 of \cite {G}. The proof is the same as the one of the latter
(this was pointed out to us by Ross Geoghegan) hence we just
describe it briefly. An action of a group on
a $CW$ complex is {\it rigid} if the stabilizer of any cell acts
trivially on that cell.
Observe that a group acting by simplicial automorphisms on a
simplicial complex acts rigidly on the first barycentric subdivision
of the complex.
Let $\pi:A\to C$ be a cellular map between
$CW$ complexes and let for each cell $e$ of $C$ a CW complex $F_e$
be given. Then (see e.g. Chapter 2.6 in \cite{G}) $\pi$ is a {\it
stack of CW complexes with fiber $F_e$} if it is (roughly
speaking) ``built" by induction on skeleta of $C$ so that over the
interior $\stackrel{\circ}{e}$ of a cell $e$, $\pi$ is homotopically equivalent to
the projection $F_e\times {\stackrel{\circ}{e}} \to {\stackrel{\circ}{e}}$.

\begin{prop}
\label{G582} Let $Y_1$ and $Y_2$ be two contractible rigid $G-CW$
complexes with cocompact $G$-actions such that stabilizers of all
cells are finite. Then $Conn_n^{\infty} (Y_1)$ iff $Conn_n^{\infty}
(Y_2)$.
\end{prop}
\dow
For $i=1,2$ apply the Borel Construction (see e.g. Chapter 2.6 in \cite{G})
using a $K(G,1)$-complex $X$ of finite type to get the commutative
diagram

$$
\begindc{\commdiag}[3]
\obj(10,15)[a]{$\widetilde X \times Y_i$}
\obj(40,15)[b]{$Y_i$}
\obj(10,1)[c]{$Z_i$}
\obj(40,1)[d]{$\Gamma _i=G\backslash Y_i$}
\mor{a}{b}{${\rm projection}$}
\mor{c}{d}{$q_i$}
\mor{a}{c}{}
\mor{b}{d}{}
\enddc
$$

Here, the diagonal action of $G$ on the contractible $CW$complex
$\widetilde X \times Y_i$ is free, so $Z_i=G\backslash (\widetilde X \times Y_i)$
is a $K(G,1)$-complex.
Then $q_i:Z_i\rightarrow \Gamma _i$ is a stack of $CW$ complexes
which can be rebuilt to give a commutative diagram

$$
\begindc{\commdiag}[3]
\obj(10,15)[a]{$Z_i'$}
\obj(40,15)[b]{$Z_i$}
\obj(25,1)[c]{$\Gamma _i$}
\mor{a}{b}{$h_i$}
\mor{b}{c}{$q_i$}
\mor{a}{c}{$q_i'$}
\enddc
$$

\noindent (see 2.6.4 in \cite{G}) in which $q_i':Z_i'\rightarrow \Gamma _i$
is a stack of $CW$ complexes and $h_i$ is a homotopy equivalence.
Here the fiber $F_e$ of $q_i$ over the cell $e$ of $\Gamma _i$ is
a $K(G_e,1)$-complex where $G_e$ is finite, and the fiber $F_e'$
of $q_i'$ over $e$ is a $K(G_e,1)$-complex of finite type (see
2.11.5 of \cite{G}). From Chapters 2.6 and 2.12 of \cite{G} one
sees that the map $\widetilde Z_i' \rightarrow \widetilde X\times Y_i
\rightarrow Y_i$ is a stack of $CW$ complexes in which the fiber
over a cell $\widetilde e$ of $Y_i$ is homeomorphic to the universal
cover of $F_e'$ (where $\widetilde e$ lies over a cell $e$ of $\Gamma
_i$). Thus that fiber is a contractible CW complex of finite type. It
follows that the indicated map $\widetilde Z_i' \rightarrow Y_i$ is a
CW-proper homotopy equivalence and hence $Conn_n^{\infty} (Y_i)$ iff
$Conn_n^{\infty} (\widetilde Z_i')$.

Thus to prove the proposition it suffices to show that
$Conn_n^{\infty} (\widetilde Z_1')$ is equivalent to $Conn_n^{\infty} (\widetilde
Z_2')$. But $Z_1'$ and $Z_2'$ are both $K(G,1)$-complexes hence
by Corollary 2.8.7 of \cite{G} they are homotopy equivalent and by
Theorem 3.1.23 of \cite{G} there exists a proper homotopy equivalence $\widetilde
Z_1'\rightarrow \widetilde Z_2'$. \kon

Hence for a group $G$ we can define a condition
$Conn_k^{\infty}(G)$ by requiring that
$Conn_k^{\infty}(X)$ holds for some, and hence for any
contractible, rigid $G-CW$ complex $X$ on which $G$ acts properly
discontinuously and cocompactly. If $Conn_k^{\infty}(G)$ for every
$k=-1,0,1,...,n$ then we say that $G$ is {\it n-connected at
infinity} (c.f. \cite{G} Chapter 5.4).

Note that, by Theorems 7.6.11 and 7.6.12 of \cite{G}, $Conn_k^{\infty}(G)$
is a quasi-isometry invariant in the sense that if $G$ and $H$ are quasi-isometric
groups then $Conn_k^{\infty}(G)$ iff $Conn_k^{\infty}(H)$.

Let $Y$ be a finite dimensional locally finite path connected
$CW$ complex, let $(L_i)_{i=1}^{\infty}$ be the filtration of $Y$
by compact subcomplexes and let $n\geq -1$ be an integer.
We denote by $Y -_c L_i$ the
largest subcomplex of $Y$ whose vertices are the vertices of
$Y\setminus L_i$.
We say that $Y$ is {\it $n$-acyclic at infinity} with respect
to a ring $R$ if the inverse system $\lk \widetilde H_k
(Y-_cL_i;R),\imath_i^j \rk$ is pro-trivial for every $-1\leq k \leq n$
(see e.g. Chapter 5.8 in \cite{G}).
Here $\imath_i^j\colon H_k(Y-_cL_j;R) \to H_k(Y-_cL_i;R)$ is a map
between reduced homologies induced by inclusion. An inverse system
like the above
is {\it pro-trivial} if for every $i$ there exists $j\geq i$
such that the image of $\imath_i^j$ is trivial.

Let moreover $Y$ be a contractible rigid $G-CW$ complex with
cocompact $G$-action such that stabilizers of all cells
are finite. If $Y$ is $n$-acyclic at infinity then we
say that the group $G$ is {\it $n$-acyclic at infinity}.
By Theorem 5.8.2 in \cite{G} this is a property of $G$,
i.e. it does not depend on the choice of $Y$ as above.

\subsection{Semistability at infinity}
\label{sei} Let $G$ act geometrically (i.e. properly discontinuously
and cocompactly) and rigidly on a one ended contractible $CW$
complex $Y$. Let $(L_i)_{i=1}^{\infty}$ be the filtration of $Y$
by compact subcomplexes. Choose a map $\omega :[0,+\infty)\to Y$ such
that $\omega (i)\in Y-_c L_i$, $i\in {\bf N}$ and put $\pi
_1^i=\pi _1(Y-_c L_i,\omega(i))$ and denote by $P_i^j:\pi
_1^{j}\rightarrow \pi_1^i$ the maps induced by inclusions and changing of
the base point along $\omega$. Then we say (see \cite{G}, Chapters
5.1, 5.4 and Theorem 5.8.2) that a group $G$ is {\it semistable at
infinity} if for every $i$ there exists $j\geq i$ such that for
all $k\geq j$ we have ${\rm im \;}(P_i^j)= {\rm im \;}(P_i^k)$.
Note that this definition does not depend on the choice of the filtration
$(L_i)_{i=1}^{\infty}$.

If $G$ is semistable at infinity then the {\it fundamental pro-group at infinity}
of $G$ being the inverse limit of the system $\lk \pi _1^i , P_i^j \rk$ is well defined
i.e. it is independent of $Y$ and $\omega$. Moreover in that case some other
invariants (e.g. {\v C}ech fundamental group or strong fundamental group)
can be defined (see Section 5.4 of \cite{G}).

It is still unknown if there exists a finitely presented one-ended group
which is not semistable at infinity.

\section{Connectedness at infinity of systolic complexes and groups}
\label{hcis}
In this section we show that for a systolic group $G$ and a
natural $k\geq 2$ the condition $Conn_k^{\infty}(G)$ holds (see
section \ref{sci} for the definition). As a consequence we get
that systolic groups are not simply connected at infinity.
Those results allow us to
distinguish systolic groups from some other natural
classes of groups.

\begin{tw}
\label{tw1}
Let $X$ be a locally finite systolic complex, let $G$ be a systolic
group and let $k\geq
2$ be a natural number. Then the conditions $Conn_k^{\infty}(X)$
and $Conn_k^{\infty}(G)$ hold.
\end{tw}
\dow
Let a compact $K\subset X$ be given. Choose a simplex $\tau$ of $X$.
We denote by $B_l$ the combinatorial ball $B_l(\tau,X)$.

There exists a natural $l>0$ such that $K$ is a subset of $B_l$.
Define $L=B_{l+2}$. Let $M=X-_cB_{l+1}$ be the maximal subcomplex
of $X$, whose vertices are the vertices of $X\setminus B_{l+1}$.
By definition $M$ is a full subcomplex of $X$ and hence,
by Fact \ref{f0.1}, it is aspherical. Therefore every map
$f:S^k\rightarrow X\setminus L\subset M$
extends to a map $g:B^{k+1}\rightarrow M \subset
X\setminus K$. Thus $Conn_k^{\infty}(X)$ holds.


To show that $Conn_k^{\infty}(G)$ holds, it is enough to notice
that if $G$ acts properly discontinuously and cocompactly by
automorphisms on a
systolic complex $Y$ then $G$ acts geometrically and rigidly
on the first barycentric
subdivision $Y'$ of $Y$.
Thus $Conn_k^{\infty}(Y)$ implies $Conn_k^{\infty}(Y')$ and
the latter is equivalent to $Conn_k^{\infty}(G)$.
\kon


\begin{tw}
\label{tw3}
Systolic groups are not simply connected at infinity.
\end{tw}
\dow
Let $G$ be a systolic group and let $X$ be a systolic simplicial
complex on which $G$ acts
geometrically, by automorphisms.

We will prove that $G$ is not simply
connected at infinity, arguing by contradiction.

Assume $G$ is simply connected at infinity. Then condition
$Conn_k^{\infty}(G)$ holds for every $k=-1,0,1,2,...$, i.e.
$G$ is $l$-connected at infinity for arbitrarily large $l$.
By the Proper Hurewicz Theorem
(Theorem 5.7.6 in \cite{G})
we have that $G$ is $l$-acyclic at infinity with respect
to ${\bf Z}$, for arbitrarily large $l$.
By Corollary 4.2 in \cite{GM1} (compare also Theorem
in \cite{GM2} and Theorem 4.3.3 in \cite{G})
it follows that $H^k(G;{\bf Z}G)=0$ for all $k$.

This is a contradiction since, by Proposition \ref{I2},
there exists a natural number $n\leq {\mr {dim}}(X)<\infty$
such that $H^n(G;{\bf Z}G)\neq 0$.
\kon

\begin{cor}
\label{co1}
Let $0\to K \to G \to H\to 0$ be a short exact sequence
of infinite finitely presented groups. If $G$ is systolic
then neither $K$ nor $H$ has one end.
\end{cor}
\dow
It follows directly from Corollary 1.5 in \cite{GM1}
(compare also Theorem 5.9.7 in \cite{G})
and the fact that $G$ is not simply connected at
infinity.
\kon

\subsection{Non-examples}
\label{nex}
Here we list some corollaries of the above results.

{\bf 1)} ${\bf Z}^3$ is not systolic (see below). This is also proved in
\cite{W} by showing that the second isoperimetric function of a
systolic group is linear (for a precise definition and a proof see
Section 8 of \cite{JS2}). Note that it is conjectured there that
higher isoperimetric functions are also linear.

{\bf 2)} Let $M$ be a closed manifold. Then $M$ has homotopy type of a
finite polyhedron $P$ (cf. \cite{KS}). The universal covers
${\widetilde M}$ and ${\widetilde P}$ of, respectively $M$ and $P$ are
properly homotopy equivalent.
Hence if ${\widetilde M}={\bf R}^n$, for
$n\geq 3$, then ${\widetilde P}$ is simply connected at infinity.
Thus we get the following.

\begin{cor}
\label{co11}
For $n\geq 3$, fundamental groups of closed manifolds covered by
${\bf R}^n$ are not systolic.
\end{cor}

{\bf 3} A theorem of D. Wise on subgroups (see \cite{W}) says that
finitely presented subgroups of a fundamental group of a compact $6$-large
complex are systolic.
Hence we get that fundamental groups of closed manifolds covered
by ${\bf R}^n$, $n\geq 3$ are not isomorphic to subgroups of
fundamental groups of compact $6$-large
complexes.

In \cite{JS2} it is proved directly (i.e. without using the result
of Wise) that fundamental groups of closed non-positively curved
riemannian manifolds of dimension at least three
are not isomorphic to subgroups of systolic
groups.



{\bf 4)} From Corollary \ref{co1} it follows that if a systolic group
$G$
is a product $G=H\times K$ of infinite groups,
then neither $K$ nor $H$ is one--ended. To show this
observe that $G$ is finitely presented
and therefore its factors $H$ and $K$ are finitely presented, too.

Corollary 7.5 in \cite{JS2} states that in fact
a product of infinite groups $H$ and $K$ is systolic only if
$H$ and $K$ are
both virtually free.

{\bf 5)} By Corollary 1.5 in \cite{GM1} (compare
Corollary 5.9.6 in \cite{G})
a product of infinite finitely generated groups has one end.
Hence, analogously to p. 4) above, we have that if
$G=H\times K \times N$, where $H$, $K$ and $N$ are infinite
then $G$ is not a systolic group.

This gives another proof of Corollary 7.7 in \cite{JS2}.


{\bf 6)} Coxeter groups acting cocompactly on $\bR$ for $n\geq 3$
are not systolic.
Moreover, let $(W,S)$ be a Coxeter system, $L$ the associated nerve
and ${\cal S}$ the poset of spherical subsets of $S$ (we use the
notation of \cite{DM}). Assume $L-\sigma$ is simply connected for
every $\sigma \in {\cal S}$.
Then by Theorem 4.3 in \cite{DM} $W$ is simply connected
at infinity and $W$ is
not systolic. (Compare point {\bf 8)} below). 

{\bf 7)} Let $X$ be a locally compact building (triangulated) with an
apartment isometric to $A$ (for basic facts about buildings see
e.g.. \cite{B}). Let $A$ be a space for which
$Conn_k^{\infty}(A)$ does not hold for some $k\geq 2$. Then
$Conn_k^{\infty}(X)$ does not hold. Hence discrete and cocompact
subgroups of automorphisms group of $X$ are not systolic.

To show this let us fix a folding $p:X\rightarrow A$ onto an
apartment of $X$ such that $p^{-1}(b)=b$ for some
chamber $b$. Take a natural $k\geq 2$ and a compact set
$\bar K\subset A$ such that for every compact
$\bar L\supset \bar K$ there exists a map
$f:S^k\rightarrow A\setminus \bar L$ non homotopically trivial in
$A\setminus \bar K$. Let $K=p^{-1}(\bar K)$ and let a compact $L
\subset X$ with $K\subset L$ be given (here we use local compactness
of $X$). Then $\bar L=p(L)$ is compact and one can find a map $f$
as above. To show that $f$ is
homotopically nontrivial in $X\setminus K$ we will argue by
a contradiction. Let us assume there exists its
extension $g:B^{k+1}\rightarrow X\setminus K$. Then
$p\circ g: B^{k+1}\rightarrow A\setminus \bar K$ is an extension of
$f:S^k\rightarrow A\setminus \bar L$ that contradict non triviality of
$f$ in $A\setminus \bar K$. In particular cocompact lattices
in groups of automorphisms of
Euclidean (i.e. with apartments being the Euclidean spaces ${\bf E}^n$)
and hyperbolic
(with apartments ${\bf H}^n$) buildings are not systolic for $n\geq 3$.

{\bf 8)} Let $X$ be a finite, locally $CAT(0)$, piecewise Euclidean
complex. For a cell $\sigma$ of $X$ let $Lk(\sigma)$ denote its
link and for $p\in Lk(\sigma)$, let $PLk(\sigma,p)$ denote the {\it
punctured link of $\sigma$ at $p$} i.e. $Lk(\sigma)$ with all the
points closer than $\pi /2$ to $p$ (in the angular metric on
$Lk(\sigma)$) removed (we use here the notation of \cite{BMcCM}).
Assume that for every cell $\sigma$ of $X$ and every point $p\in
Lk(\sigma)$ both $Lk(\sigma)$ and $PLk(\sigma,p)$ are simply
connected.
Then, by Theorem 1.3 in \cite{BMcCM}, the universal cover
$\widetilde X$ of $X$ is simply connected and hence the fundamental
group of $X$ is not systolic.
\vspace{0.2cm}

{\bf Questions.} Do there exist closed aspherical manifolds of dimension
above $2$, whose
fundamental groups are systolic ?

More generally---does there exist a closed
aspherical manifold $M$ of dimension
$n\geq 3$, such that $Conn_k^{\infty}(\widetilde M)$ holds for
every $k\geq 2$ ?

For Davis
aspherical manifolds $M$ (cf. \cite{D}),
$Conn_k^{\infty}(M)$ fails for some $k\geq 2$ (this fact we know from
Craig Guilbault and Tadeusz Januszkiewicz).


\section{Systolic chamber complexes}
\label{scis}
The aim of this section is to show that 
systolic complexes of some class are one--ended,
are not simply connected
at infinity and groups acting on them geometrically are
semistable at infinity.
That class contains, in particular, systolic groups
acting on normal pseudomanifolds of arbitrarily large
dimension,
constructed in \cite{JS1}.

Throughout this section, if it is not stated otherwise,
we denote by $X$
a systolic chamber complex of dimension $n\geq 1$.
For a natural $k$, we denote by $B_k$ the closed combinatorial ball
$B_k(\tau, X)$ around a given simplex $\tau$ of $X$.

\begin{fakt}
\label{fakt2.1}$S_k$ is an $(n-1)$-dimensional chamber
complex for every $k\geq 1$.
\end{fakt}
\dow We first show that $S_k$ is at most $(n-1)$-dimensional.
If
there existed an $l$-simplex $\sigma$ in $S_k$,
with $l \geq n$ then one would
have a projection of $\sigma$ onto $B_{k-1}$ (see Lemma \ref{l0.3})
spanning with $\sigma$ a simplex of dimension $>n$. Now we
show that every $l$-simplex $ \rho$ of $
S_k$ is contained in some $(n-1)$-simplex by induction on the dimension of
$X$. For ${ \mathrm{dim}}(X)=1$ the assertion is clear. Let now ${
\mathrm{dim}}(X)=n>1$. Consider the link $X_{ \rho}$. It is
obviously an $(n-l-1)$-dimensional chamber complex. Moreover (see
Lemma \ref{l0.3}) $X_{ \rho} \cap B_k$ is a ball $B'$ in $X_{
\rho}$ of radius $1$ around some simplex $\nu$ of $X_{ \rho}$. But
$X_{ \rho}$ is a $6$-large complex and closed balls of radii less
than $3$ in such complexes are the same as those in their
universal covers i.e. as in systolic complexes (see Lemma
\ref{l0.15}). Hence by the induction assumptions there exists
$(n-l-2)$-simplex $ \kappa$ in $
S'=
\partial(X_{\rho}\cap B_k)= X_{\rho}\cap S_k$ (see
Lemma \ref{l0.3}). Then
$\rho \ast \kappa$ is an $(n-1)$-simplex of $S_k$.

Now we show that every codimension one simplex of $S_k$
is contained in at least two maximal simplices of $S_k$.
If $\rho$ (as above) is an $(n-2)$-simplex then $X_{\rho}$ is a
$1$-dimensional chamber complex and $X_{\rho}\cap B_k=B_1(\nu,
X_{\rho})$. Since $X_{\rho}$ is a chamber complex, there should be
at least two edges in $X_{\rho}\cap B_k$ intersecting $\nu$ with
other vertices $\kappa_1$ and $\kappa_2$ lying on
$S_1(\nu, X_{\rho})$. Then $\rho \ast \kappa_1$ and $\rho \ast
\kappa_2$ are two $(n-1)$-simplices containing $\rho$. \kon
\begin{fakt}
\label{fakt2.2}Let $ \sigma$ be an $(n-1)$-dimensional simplex of
$
S_k$. Then there exists a vertex $v$ at a distance $k+1$ from
$\tau$ such that $v \ast \sigma$ is a simplex of $X$.
\end{fakt}
\dow Since $\sigma$ belongs to $S_k$, there exists a
vertex $w$ at distance $k-1$ from $\tau$ such that $w \ast
\sigma$ is a simplex. Since $X$ is a chamber complex, there is
another vertex $v$ spanning an n-simplex with $\sigma$. By Lemma
\ref{fakt2.1} the vertex $v$ does not belong to $S_k$.
Thus it
remains to show that $v$ is not at a distance $k-1$ from $\tau$.
Suppose it is. Then $P_{B_{k-1}}(b_{\sigma})$ is, by
definition, a barycenter of some simplex $\rho$ containing both $v$
and $w$. But then one has a simplex $\rho \ast \sigma$
belonging to $X$ and of dimension at least $n+1$. This
contradicts the assumption on the dimension of $X$. \kon

\begin{cor}
\label{c2.25}
The universal cover of a $6$-large chamber complex of dimension at least
$1$ is unbounded.
\end{cor}
\dow
By Lemmas \ref{fakt2.1} and \ref{fakt2.2} spheres of arbitrarily
large radii around some simplex of the universal cover
(which is systolic)
are non-empty.
\kon

Let us define here, for a given vertex $v$ of $X$ a condition
$R(v,X)$ that will be crucial for the remaining part of the section:
$$
R(v,X) \; \; {\mr {iff}} \; \;
(\forall \; \sigma \in X_v \;(X_v \setminus  {\stackrel{\circ}{B_2}}(\sigma,X_v)
{\rm \; is \; connected})).
$$
\begin{fakt}
\label{f2.3}
Suppose that condition $R(v,X)$ holds for
every vertex $v$ in $
S_k$. Then for every path $
\gamma_k=(v_1,v_2,...,v_l)$ in $(
S_k)^{(1)}$ (i.e. $v_i$'s are consecutive vertices) there exists a path

$$
\gamma_{k+1}=(w_1^1,w_1^2,...,w_1^{m(1)},z_1,w_2^1,...,z_{l-
1},w_l^1,w_l^2,...,w_l^{m(l)})
$$ in
$S_{k+1}^{(1)}$ such that:
\begin{itemize}
\item $P_{B_k}(w_i^j)$ lies in the same simplex of $S_k$ as
$v_i$, for $j=1,2,...,m(i)$, and
\item $P_{B_k}(z_i)$
lies in the same simplex of $S_k$ as the barycenter $e_i$ of the edge
$(v_i,v_{i+1})$.
\end{itemize}
\end{fakt}
\dow
By lemma
\ref{fakt2.1}, for every $i=1,2,...,l-1$ there exists an $(n-1)$-simplex
$\rho_i$ containing $(v_i,v_{i+1})$. By lemma
\ref{fakt2.2} we can find vertices $z_i$ in $
S_{k+1}$ such that $z_i
\ast
\rho_i$ span a simplex in $X$. By the definition of $P_{B_k}$ we have that
$P_{B_k}(z_i)$
lies at a distance not greater than $1$ from $e_i$.

Observe that $z_i$ belongs to $(X_{v_i}
\setminus {\stackrel{\circ}{B_2}}(
\sigma_i,X_{v_i}))\cap
(X_{v_{i+1}}
\setminus {\stackrel{\circ}{B_2}}(
\sigma_{i+1},X_{v_{i+1}}))$, where $\sigma_j=X_{v_j}\cap
B_{k-1}(\tau,X)$. Hence, by $R(v_{i+1},X)$,
one can connect $z_{i-1}$ and $z_{i}$ by
$(w_i^1,w_i^2,...,w_i^{m(i)})$. Since by definition $P_{B_k}(w_i^j)
\in B_1(v_i,
S_k)$ we conclude these projections are as desired.
\kon
\begin{fakt}
\label{f2.4}Let $S_k$ be connected and let condition
$R(v,X)$ hold for every vertex $v\in (S_k)^{(0)}$. Then
$S_{k+1}$ is connected.
\end{fakt}
\dow
Let $w,z\in (S_{k+1})^{(0)}$ and let $\mu$ and $\nu$ be the two
maximal simplices of $S_k$ containing, respectively,
$P_{B_k}(w)$ and $P_{B_k}(z)$. Choose vertices $w'$ of $\mu$
and $z'$ of $\nu$ and a path $\gamma_k$ in $(S_k)^{(1)}$
joining $w'$ and $z'$. Then by the above lemma
there exists a path $\gamma_{k+1}=(v_1,...,v_l)$ in $(S_{k+1})^{(1)}$
such that $v_1 \in X_{w'} \setminus {\stackrel{\circ}{B_2}}(\sigma_{w'}, X_{w'})$
and $v_l \in X_{z'} \setminus {\stackrel{\circ}{B_2}}(\sigma_{z'}, X_{z'})$ for appropriate $\sigma_{w'}$
and $\sigma_{z'}$ (i.e. as in Lemma \ref{l0.3}).
Since also $w\in X_{w'} \setminus {\stackrel{\circ}{B_2}}(\sigma_{w'}, X_{w'})$ and
$z\in X_{z'} \setminus {\stackrel{\circ}{B_2}}(\sigma_{z'}, X_{z'})$ by condition $R(v,X)$
we can extend $\gamma_{k+1}$ to a path connecting $w$ and $z$.
\kon

\begin{tw}
\label{tw}Let $X$ be a locally finite systolic chamber complex of
dimension greater than one. Suppose that
$R(v,X)$ holds
for almost every vertex $v\in X$ (i.e. for all but finitely many).
Then $X$ is not
simply connected at infinity. Moreover, if $G$ acts properly and
cocompactly on $X$ then $G$ is semistable at infinity.
\end{tw}
\dow
We show $X$ is not simply connected at infinity,
arguing by contradiction.
Assume $X$ is simply connected at
infinity. Let $\tau \in X^{(0)}$. Let $N\in \bN$ be such that for every
vertex $v\in X\setminus {\stackrel{\circ}{B_N}}$ condition
$R(v,X)$ holds.
By simple connectedness at infinity there exists a natural $L$ such
that every loop in $X
\setminus {\stackrel{\circ}{B_{N+L}}}$ is contractible in
$X\setminus {\stackrel{\circ}{B_N}}$. Since
$S_N$ is a full subcomplex of $X$, it is $6$-large.
It is also finite. But by Corollary \ref{c2.25}
its universal cover (i.e. a systolic complex) is infinite hence
$S_N$ is not simply connected and we can find a closed path
$\gamma_N$ in $(S_N)^{(1)}$ which is homotopically non-trivial in
$S_N$. Observe that, by Lemma
\ref{f2.3}, we can find a closed path $\gamma_{N+1}
\in (S_{N+1})^{(1)}$ such that $P_{B_N}(\gamma_{N+1})$
is homotopic to $\gamma_N$ in $S_N$.
Since $P_{B_N}$ is a deformation retraction (see \ref{l0.4}) we obtain that
$\gamma_{N+1}$ is homotopic to $\gamma_N$ within
$B_{N+1}\setminus {\stackrel{\circ}{B_N}}$. Continuing this process we can
find a closed path $\gamma_{N+L}$ in
$(S_{N+L})^{(1)}$ such that it is homotopic to $
\gamma_N$ within $B_{N+L}\setminus {\stackrel{\circ}{B_N}}$
and is homotopically trivial within $X\setminus
{\stackrel{\circ}{B_{N+L}}}$
by assumptions on $L$. But then $\gamma_N$ is homotopic in
$S_N$ to $P_{B_N}(\gamma_{N+L})$ which is
homotopically trivial within $S_N$. But this contradicts
the choice of $\gamma_N$.

To see that $G$ is not simply connected at infinity observe
that $G$ acts properly discontinuously, cocompactly and rigidly
on the barycentric subdivision of $X$ (see Section \ref{sci}).

For semistability take a filtration $(B_k)_{k=0}^{\infty}$ of $X$
(see Section \ref{sei}) and a suitable ray $\omega :[0,+\infty)\rightarrow
X$ such that $P_{B_k}(\omega (k))=\omega (k-1)$ (this can be done
by taking a sequence $(x_k)_{k=0}^{\infty}$ of points such that
$x_k\in S_k$ then considering their projections on spheres
$S_i$ which are compact and getting the desired sequence by a diagonal
argument taking accumulation points of projections on $S_i$'s).
Observe that $X'-_cB_k$ (we take a barycentric
subdivision $X'$ of $X$ to make action of $G$ rigid) is homotopically
equivalent to $X\setminus {\stackrel{\circ}{B}}_k$ and hence to $S_{k+1}$.
Given $k\geq N$ ($N$ as above)
and an element $g\in \pi _1(X-_cB_k,\omega (k))=
\pi _1(S_{k+1},\omega (k))$ choose a path $\gamma$ in
$(S_{k+1})^{(1)}$ representing $g$. Then by the above construction
for any $L=0,1,2,...$ one can find a closed path
$\gamma _L$ in $(S_{k+L+1})^{(1)}$ such that
$P_{B_{k+1}}(\gamma _L)$ is homotopic to $\gamma$ in $S_{k+1}$.
Hence if $g_L$ is an element of $\pi _1(S_{k+L+1},
\omega (k+L))=\pi _1(X-_cB_{k+L},\omega (k+L))$
represented by $\gamma _L$ (using connectedness
of $S_{k+L+1}$ by Lemma \ref{f2.4}) we have
that $(P_{B_{k+1}})_{\ast}(g_L)=g$ and that the map
$\pi _1(X-_cB_{k+L},\omega (k+L))\rightarrow \pi _1(X-_cB_k,\omega (k))$
induced by the inclusion is surjective.
This shows $G$ is semistable at infinity.
\kon

The next lemma gives a condition which helps to prove the condition
$R(v,X)$ in some cases---see e.g. Section \ref{extra}.

\begin{fakt}
\label{f2.5}Let
$X$ be a 6-large chamber complex  such
that
$X_{\kappa}$ is connected for every simplex $\kappa$ of
codimension greater than one in $X$
and
$X_{\sigma}\setminus {\stackrel{\circ}{B_2}}(\rho,X_{\sigma})$
is connected for every codimension two simplex $\sigma$ of $X$ and every simplex
$\rho$ of its link $X_{\sigma}$. Then $R(v,X)$ holds
for every vertex $v$ of $X$.
\end{fakt}
\dow
We will proceed by induction on $n={\rm dim}(X)$.

For $n=2$ the assertion is clear since codimension two simplices
are just vertices.

Assume we proved the lemma for $n\leq k$. Let ${\rm dim}(X)=k+1$.
Take a vertex $v$ of $X$ and consider its link $X_v$. It has
dimension $k$. Observe that if
$\alpha$ is a codimension $l$ simplex of $X_v$ then
$\beta=\alpha\ast v$ is a codimension $l$ simplex of $X$ and
$(X_v)_{\alpha}=X_{\beta}$.
Hence by the induction assumptions for every vertex $w$ of $X_v$
condition
$R(w,X)$ holds.

Take a simplex $\omega$ of $X_v$. Since $\omega=
S_0(\omega,X_v)$ is connected we have by Lemma \ref{f2.4},
Lemma \ref{l0.15} and by
the fact that $X_v$ is $6$-large (as a full subcomplex of $X$)
that also $S_1(\omega,X_v)$ and
$S_2(\omega,X_v)$ are connected. Take two vertices $t$ and $s$ of
$X_v \setminus {\stackrel{\circ}{B_2}}(\omega,X_v)$. Since, by assumptions $X_v$
is connected there exists a path in $(X_v)^{(1)}$ joining them. If
this path misses $S_2(\omega,X_v)$ it joins these
vertices in $X_v\setminus {\stackrel{\circ}{B_2}}(\omega,X_v)$. If not we can
replace it (using connectedness of $S_2(\omega,X_v)$) by a
path in
$X_v\setminus {\stackrel{\circ}{B_2}}(\omega,X_v)$.
Hence we have the conclusion.
\kon
\begin{cor}
\label{c1}Let $X$ be locally finite systolic chamber complex of
dimension above one. If $X_{\kappa}$ is connected
for every simplex $\kappa$ of codimension greater than one in $X$
and $X_{\sigma}\setminus {\stackrel{\circ}{B_2}}(\tau,X_{\sigma})$ is connected
for every codimension two simplex $\sigma$ of $X$ and every
simplex $\tau$ of its link $X_{\sigma}$ then $X$ is not simply
connected at infinity and every group acting properly and
cocompactly on $X$ is semistable at infinity.
\end{cor}
\re Observe that it follows from the proof of Lemma
\ref{f2.5} that it is enough to have the above assumptions only for almost
every codimension two simplex $\sigma$ similarly as in Theorem
\ref{tw}
\begin{cor}
\label{c2}
Locally finite normal systolic pseudomanifolds are one-ended
and are not simply
connected at infinity. Groups acting on them cocompactly and
properly are one--ended, semistable at infinity and
do not split as amalgamated products or as HNN extensions over
finite groups.
\end{cor}
\dow
Let $X$ be a locally finite normal systolic pseudomanifold.
Let $\sigma$ be an $l$-simplex of $X$, for $l\geq 2$.
Then $\partial \sigma$ is connected and by Lemma \ref{f2.4}
every sphere $S_i(\sigma, X)$ is connected. Thus
$X$ has one end.

As for semistability at infinity
observe that
one dimensional link of a normal pseudomanifold is a circle.
Hence it satisfies the assumptions of the preceding corollary.

By Stallings theorem (\cite{S}) if a finitely generated group has one
end then it does not split as an amalgamated product or an HNN extension over
a finite group.
\kon

{\sl Remark.} Systolic groups acting geometrically on normal
pseudomanifolds constructed in \cite{JS1} are the only systolic
groups of cohomological dimension greater then two known to us
at the moment.

\vspace{0.3cm}
{\bf Example.}
We give an example of a $2$-dimensional systolic normal
chamber complex $X$ which is simply connected at infinity.
It shows we cannot delete condition $R(v,X)$ from
the assumptions of Theorem \ref{tw}.
Let $(A,a,b,r)$ denote the quadruple consisting of:

- $A$ - a closed Euclidean half-plane (i.e. $\left\{ (x,y) \in {\bf
R}^2:y\geq0\right\}$) triangulated by regular triangles,

- $a$ - a vertex (of triangulation) in
the interior of $A$,

- $b$ - a vertex on the boundary $\partial A$,

- $r$ - a closed geodesic ray in $A^{(1)}$ starting from $a$.

Let $(A_i,a_i,b_i,r_i)_{i=1}^{\infty}$ be a sequence of quadruples
isomorphic to $(A,a,b,r)$. We will construct $X$ recursively by
``gluing consecutive quadruples". More precisely let $X_1=A_1$.
Having $X_k$ we find the lowest index $i_k\leq k$ such that
$b_{i_k}$ belongs to a boundary edge $e_k$ of $A_{i_k}$ (viewed as
a subspace of $X_k$) that does not belong to any triangles of
$A_j$ (viewed as a subspace of $X_k$), for $j\neq i_k$, $j\leq k$.
Then $X_{k+1}=X_k\cup_{\psi_k} A_{k+1}$, where
$\psi_k:r_{k+1}\rightarrow A_{i_k}$ sends isometrically $r_{k+1}$
to a closed half-line $\bar r_{k+1}$ of $\partial A_{i_k}$
starting at $b_{i_k}$ and containing $e_k$. If one defines a map
$\phi_k:X_{k+1}\rightarrow X_k$ so that $\phi_k(A_{k+1})= \bar
r_{k+1}$ and as identity on $X_k\subset X_{k+1}$ then
$(X_i,\phi_i)_{i=1}^{\infty}$ is an inverse sequence and we set
$X={ \mathrm {inv \; lim}}(X_i,\phi_i)$.

\vspace{0.3cm}
\epsfxsize=11cm
\frame{\epsffile{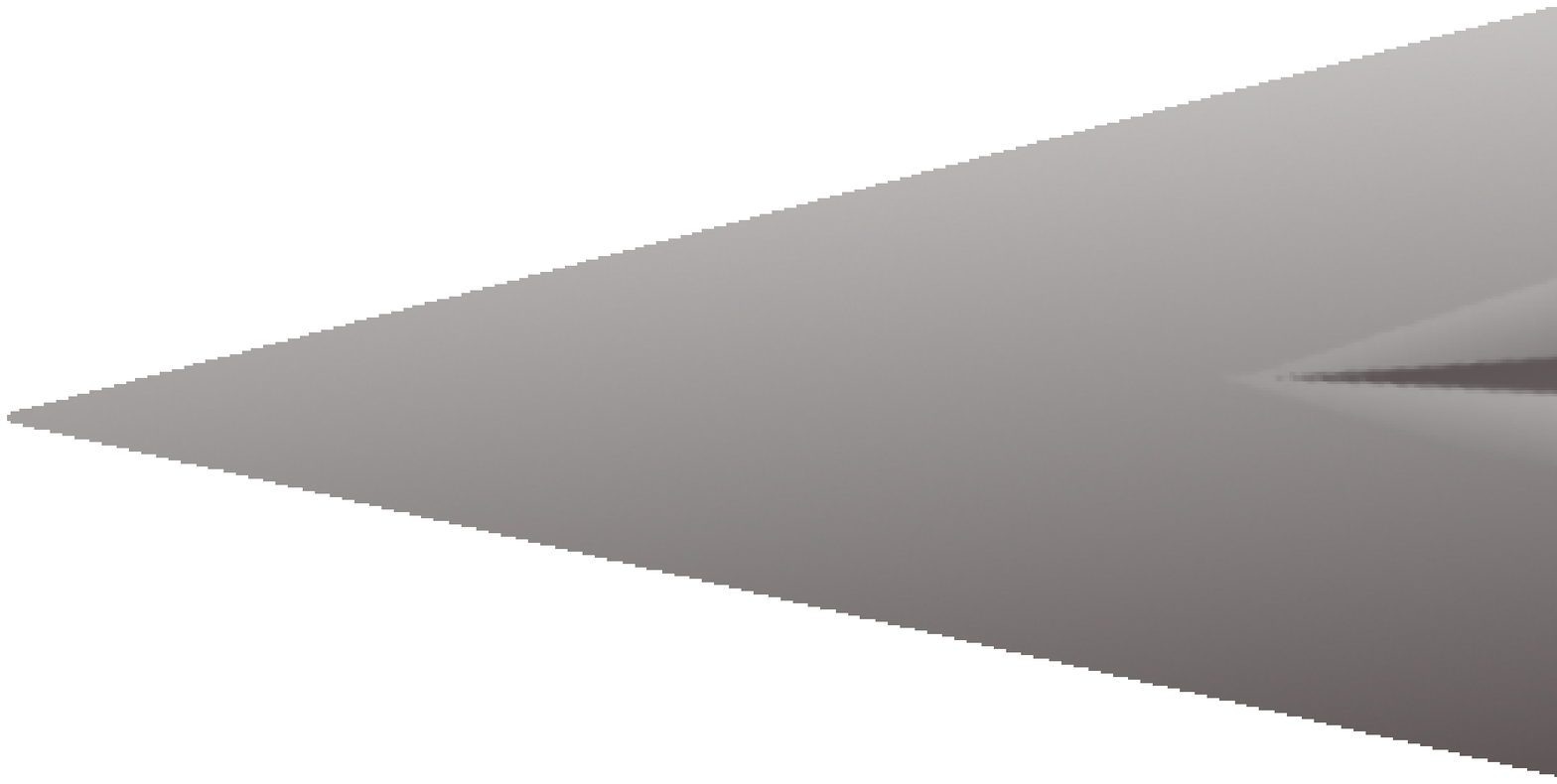}}
\begin{center} Topology of $X_{11}$ \end{center}
\vspace{0.3cm}

It is clear that $X$ is a locally finite $2$-
dimensional normal systolic chamber complex. Since every $A_i$ is contractible
and we ``glue them along rays" we obtain that $X$ is contractible. Then it can be
easily observed that $X$ is simply connected at infinity.
Observe that the assumptions of
Theorem
\ref{tw} (i.e. the condition $R(b_i,X)$)
are not satisfied for vertices $b_i$ in $X$.
Note that there is no geometric
group action on $X$.

{\bf Question.}
Are one-ended
systolic groups semistable
at infinity ?

\section{Extra-tileable simplices of groups}
\label{extra}
Besides groups acting geometrically on trees and on 2-dimensional systolic
complexes the only other systolic groups known to us at the moment are
the ones constructed in \cite{JS1} as fundamental groups of locally
$6$-large and locally extra-tilable simplices of groups (see the comment
after Proposition \ref{p0.6} in Section
\ref{pre}).

In this section we show the latter groups are one--ended and
semistable at infinity,
provided groups associated with codimension one faces are nontrivial.
Let ${\cal G}$ be a locally
$6$-large and locally extra-tilable simplex of groups, as constructed
in Proposition 19.1 of \cite{JS1}
and let $\Delta$ be the underlying simplex.
Assume that groups associated with codimension one faces
of $\Delta$ are nontrivial.
In view of Corollary \ref{c1} it is enough to show that
the universal development $\widetilde {\cal G}$ of ${\cal G}$
(which exists by Proposition \ref{p0.5})
is a locally finite systolic chamber complex such that its links
of dimension above $0$ are connected and complements of open
balls of radius $2$ in $1$-dimensional links of $\widetilde {\cal G}$
are connected.
The argument below was supplied by
J. {\' S}wi{\c a}tkowski.

We refer to the construction made in \cite{JS1}, Section 19, and
use notation introduced there. See also Section \ref{pre} above.

By Proposition 17.1 (5) of \cite{JS1} (compare also \cite{BH})
the universal development
$\widetilde {\cal G}$ is a union of its top-dimensional simplices
and, by our assumption,
the groups $G_{\sigma}$ (compare Section \ref{pre})
are non-trivial, for all codimension one faces $\sigma$ of
$\Delta$. Hence $\widetilde {\cal G}$ is a chamber complex.
By Proposition 19.1 of \cite{JS1} it is locally $6$-large
and hence---because it is simply connected---systolic.
It is locally finite by Proposition 19.1 (3) of \cite{JS1},
because the groups $G_{\sigma}$ are finite, for all faces
$\sigma$ of $\Delta$ (by the proof of
Proposition 19.1 of \cite{JS1}). By construction
(the proof of
Proposition 19.1 of \cite{JS1}), for every face $\sigma$ of
$\Delta$ the canonical morphism $i_{\sigma}:{\cal G}^{\sigma}
\to G_{\sigma}$ is surjective (meaning $G_{\sigma}$ is generated
by the union of images of $i_{\sigma}$'s). Hence, by Proposition 17.1
(6) of \cite{JS1}, links in $\widetilde {\cal G}$ are gallery connected.

Now, we show that complements of open
balls of radius $2$ in $1$-dimensional links of $\widetilde {\cal G}$
are connected. Let $\sigma$ be a codimension $2$ face of $\Delta$
and $s$, $t$ be the two codimension $1$ faces containing $\sigma$.
Then the groups $G_{s}$,
$G_t$ and $G_{\Delta}$, form a segment of groups ${\cal
G}^{\sigma}$ over the link $\Delta_{\sigma}$ which is
extra-tilable, by local extra--tilability of ${\cal G}$.
The link of a codimension two simplex of the universal
development $\widetilde {\cal G}$ of ${\cal G}$ lying over $\sigma$ is
constructed as follows (cf. Corollary 18.6 in \cite{JS1}). In the
tree $\widetilde {{\cal G}^{\sigma}}$ (the universal development of
${\cal G}^{\sigma}$) take a ball $B=B_k([\Delta_{\sigma},1],\widetilde
{{\cal G}^{\sigma}})$ of radius $k=6$. Let $H_{B}$ be the subgroup
of the direct limit $\widehat {{\cal G}^{\sigma}}$ of ${\cal G}^{\sigma}$
for which $B$ is a strict fundamental domain.
Such a subgroup exists by Proposition \ref{p0.6}.
The $1$-dimensional link we look for,
is of the form $N\backslash \widetilde {{\cal G}^{\sigma}}$,
for some finite index
normal subgroup
$N$ of $H(B)$.
The image $p(B)\subset N\backslash \widetilde
{{\cal G}^{\sigma}}$ of $B$ is a strict fundamental domain for the action
of a finite group $H_B/N$ on $N\backslash \widetilde {{\cal G}^{\sigma}}$,
where $p:\widetilde {{\cal
G}^{\sigma}}\to N\backslash \widetilde {{\cal G}^{\sigma}}$
is the quotient map.

We show that
$A:=(N\backslash \widetilde {{\cal G}^{\sigma}}) \setminus p(
{\stackrel{\circ}{B_2}([\Delta_{\sigma},1],\widetilde {{\cal G}^{\sigma}})})$ is connected.
Then $(N\backslash \widetilde {{\cal G}^{\sigma}}) \setminus p(
{\stackrel{\circ}{B_2}([\Delta_{\sigma},g],\widetilde {{\cal G}^{\sigma}})})$
is connected for every $g\in \widehat {{\cal G}^{\sigma}}$ since
$N\backslash \widehat {{\cal G}^{\sigma}}$ acts transitively on
the set of top-simplices of the link $N\backslash \widetilde {{\cal G}^{\sigma}}$.
We show that every two points $x,y\in A$ can be connected by a path
in $A$.
Consider first the case when $x\in A\setminus p(B)$ and $y\in A\cap p(B)$.
Consider the Cayley graph $(H_B/N,S)$ where $S$ is the set of generators
defined as follows
$$
s\in S \; \; {\rm iff} \; \; p(B)\cap s(p(B))\neq \emptyset.
$$
This set generates $H_B/N$ due to connectedness of $B$ and
$N\backslash \widetilde{{\cal G}^{\sigma}}$.
Find a path in $A\cap p(B)$ joining $y$ with $y_0\in p(\partial B)$.
Let $1\neq s\in T$ be such that $y_0\in s(p(B))$.
It is a general easy fact that the complement of a vertex
(in particular, the complement of the vertex corresponding to $1$) in any
Cayley graph of a finite group is connected. Hence if $k\in H_B/N$
is such that $x\in k(p(B))$ we can find a sequence $s,s_1,...,
s_l\in S$ such that $ss_1...s_l=k$ and $ss_1...s_i\neq 1$, for
$i=1,2,...,l$. Then $ss_1...s_i(p(B))\subset A$, for $i=1,2,...,l$.
For $i=1,2,...,l-1$ and $y_i=ss_1...s_{i-1}(p(B))\cap
ss_1...s_{i}(p(B))$ we can find a path in $ss_1...s_{i-1}(p(B))$
connecting $y_{i-1}$ and $y_i$ for $i<l$ and a path
in $k(p(B))$ connecting $y_{l-1}$ and $x$. Concatenation of those paths
is a path in $A$ connecting $x$ and $y$. If $x,y\in A\cap p(B)$,
then we can find a vertex $z\in A\setminus p(B)$ and connect
$x$ and $y$ with $z$.

Similarly for $x,y\in A\setminus p(B)$.
Thus $A$ is connected.

The proof that $\widetilde {{\cal G}^{\sigma}}/N \setminus p(
{\stackrel{\circ}{B_2}([z,1],\widetilde {{\cal G}^{\sigma}})})$ is connected
for a vertex $z$ of $\Delta_{\sigma}$ is similar.


\end{document}